%% file: dual-control.tex
\documentclass{article}

\include{dual-control-macros}

\begin{document}

\title{On Coarse Spectral Geometry in Even Dimension}
\author{Robert Yuncken}

\maketitle

\begin{abstract}
Let $\sigma$ be the involution of the Roe algebra $\Roe{\RR}$ which is induced from the reflection $\RR\to\RR; ~x\mapsto -x$.  A graded Fredholm module over a separable $C^*$-algebra $A$ gives rise to a homomorphism $\tilde{\rho}:A\to\Roe{\RR}^\sigma$ to the fixed-point subalgebra.  We use this observation to give an even-dimensional analogue of a result of Roe.  Namely, we show that the $K$-theory of this symmetric Roe algebra is $K_0(\Roe{\RR}^\sigma)\cong\ZZ$, $K_1(\Roe{\RR})=0$, and that the induced map $\tilde{\rho}_*:K_0(A) \to \ZZ$ on $K$-theory gives the index pairing of $K$-homology with $K$-theory.
\end{abstract}


\section{Introduction}

In \cite{Roe:dual-control}, Roe observed that a Dirac operator $D$ on an odd-dimensional closed manifold $M$ gives rise to a $C^*$-algebra homomorphism
\begin{equation}
  \tilde{\rho}:C(M) \to C^*|\RR|
\end{equation}
from the continuous functions on $M$ to the Roe algebra of the real line $\RR$.  The space $\RR$ appears because, up to coarse equivalence, it is the spectrum of the self-adjoint operator $D$.
The $K$-theory of $\Roe{\RR}$ is
$$
  K_n(C^*|\RR|) \cong \begin{cases} 0, & n=0,\\ \ZZ, & n=1, \end{cases}
$$
and the map 
\begin{equation}
  \tilde{\rho}_* : K_1(C(M)) \to K_1(C^*|\RR|) \cong \ZZ
\end{equation}
agrees with the index pairing of $K$-theory with the $K$-homology class $[D] \in K_1(M)$.

This point of view was extensively developed by Luu (\cite{Luu:thesis}), who showed that analytic $K$-homology can be reformulated entirely in the language of coarse spectral geometry.  Specifically, let $A$ be a separable
$C^*$-algebra.  Luu defined groups $KC^n(A,\CC)$ whose cycles are $*$-homomorphisms $\rho:A \to C^*|\RR^n|$,\footnote
{
  The most natural coarse structure on $\RR^n$ here is the topologically controlled coarse structure associated to the compactification of $\RR^n$ by a sphere at infinity.  (See \cite{Roe:coarse-geometry} for the definition.)  If $A$ is separable, it turns out to be equivalent to use the standard metric coarse structure on $\RR^n$, although the construction becomes somewhat more technical.  The $K$-theory of $C^*|\RR^n|$ is the same in either case.  
}
and then proved that $KC^n(A,\CC) \cong KK^n(A,\CC)$.  In fact, Luu worked with an arbitrary ($\sigma$-unital) coefficient algebra $B$, to produce groups $KC^n(A,B)$ isomorphic to $KK^n(A,B)$.  We choose not to work in that generality here.

Luu's picture of $K$-homology is aesthetically very pleasing.  The price of this elegance, however, is some computational complexity in even dimensions.  The isomorphism of $KK$ and $KC$ in even dimension is achieved via a map $KK^0(A,\CC) \to KC^2(A,\CC)$ which requires as input a {\em balanced} Fredholm module, \emph{i.e.}~a graded Fredholm module of the form $(H=H_0\oplus H_0, \rho=\rho_0\oplus\phi_0, F= \scriptstyle \left(\begin{smallmatrix}0&U^*\\U&0 \end{smallmatrix}\right)$ for some Hilbert space $H_0$, representation $\rho_0$ and Fredholm operator $U:H_0\to H_0$.  While every $K^0$-class can be represented by a balanced Fredholm module, the process of ``balancing'' is quite heavy-handed.  For instance, given a Dirac operator on an even dimensional manifold, the Hilbert space of the associated balanced Fredholm module is an infinite direct sum of $L^2$-sections of the spinor bundle.  (See  \cite[Proposition 8.3.12]{HR:book}.)  The relationship between the spectrum of $U$ and that of the original operator $D$ is not obvious.

In this paper, we describe an alternative approach to controlled spectral geometry in even dimension which is more convenient for geometric applications.  Let $s:\RR\to\RR$ denote the reflection through the origin.  This induces a $*$-involution $\sigma$ of the Roe algebra $C^*|\RR|$ (see Section \ref{sec:symmetric-Roe-algebra}).   
Given a graded Fredholm module $(H,\rho,D)$ for $A$, Roe's construction in fact produces a $*$-homomorphism $\tilde{\rho}:A\to C^*|\RR|^\sigma$ into the fixed-point algebra of $\sigma$.  Our main result is the following.

\begin{theorem}
\label{thm:main-theorem}
The $K$-theory of the symmetric Roe algebra is
$$
  K_n(C^*|\RR|^\sigma) \cong \begin{cases} \ZZ, & n=0,\\ 0, & n=1, \end{cases}
$$
and the induced map
\begin{equation}
  \tilde{\rho}_* : K_0(A) \to K_0(C^*|\RR|^\sigma) \cong \ZZ
\end{equation}
agrees with the index pairing of $[(H,\rho,D)]\in K^0(A)$ with $K$-theory.
\end{theorem}

\medskip

The author would like to thank Vi\^et-Trung Luu for stimulating chats.


\section{Preliminaries: The Roe algebra $C^*|\RR|$}

We shall use $|\RR|$ to denote the real line equipped with the topological coarse structure induced from the two-point compactification $\overline{\RR} \defeq \RR\cup\{\pm\infty\}$.  Thus, a set $E\subseteq\RR\times\RR$ is {\em controlled} if for any sequence $(x_n,y_n)\in E$, $x_n \to \infty$ (resp.~$-\infty$) if and only if $y_n \to \infty$ (resp.~$-\infty$).

We shall refer to a Hilbert space $H$ equipped with a nondegenerate representation $m:C_0(\RR)\to\scrB(H)$ as a {\em geometric $\RR$-Hilbert space}.  By the spectral theorem, $m$ extends naturally to the algebra of Borel functions $B(\RR)$.  We shall typically suppress mention of $m$ in the notation.  We use $\chi_Y$ to denote the characteristic function of a subset $Y\subset\RR$.

An operator $T\in\scrB(H)$ is {\em locally compact} if $fT,~Tf \in \scrK(H)$ for all $f\in C_0(\RR)$.  It is {\em controlled} (for the above topological coarse structure) if for all $R\in \RR$ there exists $S\in\RR$ such that
$$
  \begin{array}{rclcrcl}
  \chi_{(-\infty,R]} T \chi_{[S,\infty)} &=& 0,  & \quad &
    \chi_{[S,\infty)} T \chi_{(-\infty,R]} &=& 0,  \\
  \chi_{(-R,\infty]} T \chi_{(-\infty,-S]} &=& 0,  &&
    \chi_{(-\infty,-S]} T \chi_{[-R,\infty)} &=& 0.  \\
  \end{array}
$$
One defines $C^*(|\RR|;H)$ as the norm-closure of the locally compact and controlled operators on $H$.
This $C^*$-algebra is independent of the choice of $H$ as long as $H$ is {\em ample}, {\em i.e.} $m(f)$ is noncompact for all nonzero $f\in C_0(\RR)$.  In that case, the algebra is referred to as the {\em Roe algebra} $\Roe{\RR}$.

The following standard facts are easy consequences of the definitions.  The reader familiar with Roe algebras may prefer to recognize them as consequences of the coarsely excisive decomposition $\RR = (-\infty,0]\cup[0,\infty)$, where we note that the ideal $C^*_{|\RR|}(|\{0\}|;H)$ associated to the inclusion of a point into $\RR$ is just the compact operators. (See \cite{HRY},\cite{HPR}.)

\begin{lemma}
\label{lem:coarse-decomposition}
  Let $T \in C^*(|\RR|;H)$.  For any $R_1, R_2\in \RR$,
  \begin{enumerate}
    \item $\chi_{(-\infty,R_1]}T\chi_{[R_2,\infty)}$ and $\chi_{[R_2,\infty)}T\chi_{(-\infty,R_1]}$ are compact operators.
    \item $[T, \chi_{(-\infty,R_1]}]$ and $[T, \chi_{[R_2,\infty)}]$ are compact operators.
  \end{enumerate}
\end{lemma}


\section{Graded Fredholm modules and the symmetric Roe algebra}
\label{sec:symmetric-Roe-algebra}

In what follows, we shall use the unbounded (`Baaj-Julg') picture of $K$-homology.  This is a purely aesthetic choice---see Remark \ref{rmk:bounded-picture} for the construction using bounded Fredholm modules.

Let $A$ be a $C^*$-algebra, and let $(H,\rho,D)$ be a graded unbounded Fredholm module for $A$, \emph{i.e.}~$H$ is a $\ZZ/2\ZZ$-graded Hilbert space, $\rho$ is a representation of $A$ by even operators on $H$, and $D$ is an odd self-adjoint unbounded operator on $H$ such that
\begin{enumerate}
\item[(1)] for all $a\in A$, $(1+D^2)^{-\half} \rho(a)$ extends to a compact operator,
\item[(2)] for a dense set of $a\in A$, $[D,\rho(a)]$ is densely defined and extends to a bounded
operator.
\end{enumerate}
Let $\gamma_{\ev}, \gamma_\od$ denote the projections onto the even and odd components of $H$, and $\gamma=\gamma_\ev - \gamma_\od$ be the grading operator.  Let $\sigma$ be the involution of $\scrB(H)$ defined by $\sigma:T\mapsto \gamma T\gamma$.

Functional calculus on the operator $D$ provides $H$ with a geometric $\RR$ structure, namely $m: B(\RR) \to \scrB(H);~ f \mapsto f(D).$    For any $f\in C_0(\RR)$, 
$$
  \sigma(m(f)) = f(\gamma.D.\gamma) = f(-D) = m(f\circ s),
$$
where $s:\RR \to \RR$ is the reflection in the origin.  In coarse language, $\gamma$ is a covering isometry for $s$.  It follows that $\sigma$ restricts to an involution of $C^*(|\RR|;H)$.  The subalgebra fixed by $\sigma$ will be denoted $C^*(|\RR|;H)^\sigma$.

Taking this symmetry into account gives an immediate strengthening of Roe's construction for ungraded Fredholm modules.

\begin{proposition}
\label{prop:A-is-controlled}
The image of $\rho$ lies in $C^*(|\RR|;H)^\sigma$.
\end{proposition}

\begin{proof}
The function $f(x) = x(1+x^2)^{-\half}$ generates $C(\overline{\RR})$, and the ideal generated by $g(x) = (1+x^2)^{-\half}$ is $C_0(\RR)$.  Using \cite[Theorem 6.5.1]{HR:book}, Properties (1) and (2) above imply that $\rho(a)\in C^*(|\RR|;H)$ for any $a\in A$.  Since $\rho(a)$ is even, $\sigma(a) = \gamma \rho(a) \gamma = \rho(a)$.
\end{proof}

This geometric $\RR$-Hilbert space $H$ is not typically ample.  
However, one can always embed $H$ into an ample geometric $\RR$-Hilbert space.  For specificity, let us put $\scrH\defeq H \oplus L^2(\RR)$, where $L^2(\RR)$ has its natural geometric $\RR$-structure.  Extension of operators by zero gives an inclusion $\iota:C^*(|\RR|;H) \hookrightarrow \Roe{\RR}$.  Put $\tilde{\rho} = \iota \circ \rho:A \to \Roe{\RR}$.

The symmetry $g\mapsto g\circ s$ defines a grading operator on $L^2(\RR)$.  We shall reuse $\gamma$ to denote the total grading operator on $\scrH$.  Likewise, we use $\sigma$ to denote conjugation by $\gamma$ in $\scrB(\scrH)$.  Then $\tilde{\rho}$ has image in $\Roe{\RR}^\sigma$.

\begin{remark}
\label{rmk:canonicity}

In the above, we have employed a specific choice of symmetry $\sigma \in \Aut(\Roe{\RR})$ associated to the reflection $s$ of $\RR$.  For the expert concerned about the uniqueness of this definition, we supply some brief comments without proof.  They shall not be needed in what follows. 

Let $\scrH$ be any ample geometric $\RR$-Hilbert space.  By \cite[Prop.~2.2.11(iii)]{Luu:thesis} (following \cite{HRY}), there exists a unitary $\gamma:\scrH\to\scrH$ which covers $s$, in the sense that $(1\times s)(\Supp(\gamma)) \subseteq \RR\times\RR$ is a controlled set.  By carrying out the proof of this fact in a way that maintains the reflective symmetry, one can ensure that $\gamma$ is involutive, $\gamma^2=1$.  Then $\sigma:T \to \gamma T \gamma$ is an involution of $\Roe{\RR}$.  If $\gamma'$ is another involutive covering isometry for $s$, then there is a controlled unitary $V\in \scrB(\scrH)$ such that $\gamma' = V\gamma$ (\cite[Prop.~2.2.11(iv)]{Luu:thesis} following \cite{HRY}).  If $\sigma'$ is conjugation by $\gamma'$, then $\Roe{\RR}^{\sigma'} = V \Roe{\RR}^\sigma V^*$.  Thus the symmetric Roe algebra $\Roe{\RR}^\sigma$ is unique up to controlled unitary equivalence.

\end{remark}

\begin{remark}
\label{rmk:bounded-picture}
The bounded Fredholm module corresponding to $(H,\rho,D)$ is\linebreak $(H,\rho,F\defeq D(1+D^2)^{-\half})$.  The map $\phi:x\mapsto x(1+x^2)^{-\half}$ defines a coarse equivalence from $|\RR|$ to the interval $|(-1,1)|$, with topological coarse structure associated to its two-point compacification $[-1,1]$.  Thus, the bounded picture of $K$-homology provides a morphism $\rho:A\to \Roe{(1,-1)} \cong \Roe{\RR}$.
\end{remark}


\section{$K$-theory of the symmetric Roe algebra $\Roe{\RR}^\sigma$}

\begin{proposition}
\label{prop:K-theory}
The $K$-theory of $\Roe{\RR}^\sigma$ is
$$
  K_\bullet(C^*|\RR|^\sigma) \cong
    \begin{cases}
      \ZZ, & \bullet = 0,\\
      0, & \bullet = 1.
    \end{cases}
$$
Moreover, $K_0(\Roe{\RR})^\sigma$ is generated by finite rank projections $p \in M_n(C^*|\RR|^\sigma)$, and for such projections, the map to $\ZZ$ is given by
$$
  [p] \mapsto \dim p\scrH_\ev - \dim p\scrH_\od.
$$
\end{proposition}

We use a Mayer-Vietoris type argument (\emph{cf.}~\cite{HRY}).  Put $Y_+\defeq [1,\infty)$, $Y_-\defeq (-\infty,-1]$, with their coarse structures inherited from $|\RR|$.  We will abbreviate $\chi_{Y_\pm}$ as $\chi_\pm$.  Since $\scrH_+\defeq \chi_+\scrH$ is an ample geometric $Y_+$-Hilbert space, we can define the Roe algebra $\Roe{Y_+}$ as the corner algebra $C^*(|Y_+|;\scrH_+) = \chi_+\Roe{\RR}\chi_+$.  Likewise for $\Roe{Y_-}$.  

Note that $\sigma(\chi_\pm) = \chi_\mp$, so that $\sigma$ interchanges $\Roe{Y_+}$ and $\Roe{Y_-}$.  Since $\chi_+\chi_- =0$, the symmetrization map $\average:T \mapsto T+ \sigma(T)$ is a $*$-homomorphism from $\Roe{Y_+}$ into $\Roe{\RR}^\sigma$.  We obtain a morphism of short-exact sequences,
\begin{equation}
\label{eq:short-exact-sequences}
  \xymatrix{
  0 \ar[r] &
  \scrK(\scrH_+) \ar[r]  \ar[d]_{\average} &
  \Roe{Y_+} \ar[r]  \ar[d]_{\average} &
  \Roe{Y_+}/ \scrK(\scrH_+) \ar[r]  \ar[d]_{\average} &
  0 \\
  0 \ar[r] &
  \scrK(\scrH)^\sigma \ar[r] &
  \Roe{\RR}^\sigma \ar[r] &
  \Roe{\RR}^\sigma / \scrK(\scrH)^\sigma \ar[r] &
  0.
  }
\end{equation}

\begin{lemma}
  The right-hand map $\average: \Roe{Y_+}/ \scrK(\scrH_+) \to  \Roe{\RR}^\sigma / \scrK(\scrH)^\sigma$ is an isomorphism.
\end{lemma}

\begin{proof}
Let $\psi:\Roe{\RR}^\sigma \to \Roe{Y_+}$ denote the cut-down map $T \mapsto \chi_+ T \chi_+$.  By using Lemma \ref{lem:coarse-decomposition}(ii),  $\psi$ is a homomorphism modulo compacts, so it descends to a homomorphism $\psi:\Roe{\RR}^\sigma / \scrK(\scrH)^\sigma \to \Roe{Y_+}/ \scrK(\scrH)$.  By Lemma \ref{lem:coarse-decomposition}(i), for any $T\in \Roe{\RR}^\sigma$ we have
\begin{eqnarray*}
  T &\equiv& \chi_+ T \chi_+ + \chi_- T \chi_- \mod{\scrK(\scrH)^\sigma}, \\
\end{eqnarray*}
so that $\psi$ is inverse to $\average$.
\end{proof}

Put $\scrH_{\ev} \defeq \gamma_{\ev} \scrH$, $\scrH_{\od} \defeq \gamma_{\od} \scrH$.

\begin{lemma}
\label{lem:K-theory-of-KHsigma}
  We have $\scrK(\scrH)^\sigma \cong \scrK(\scrH_\ev) \oplus \scrK(\scrH_\od)$ via   $T \mapsto T\gamma_\ev \oplus T\gamma_\od$.  In particular, $K_0(\scrK(\scrH)^\sigma) \cong \ZZ\oplus\ZZ$ via the map which sends the class of a projection $p$ to $(\dim(p\scrH_\ev), \dim(p\scrH_\od))$.
\end{lemma}

\begin{proof}
Note that any $T\in\Roe{\RR}^\sigma$ commutes with $\gamma$, so $T\mapsto T\gamma_\ev\oplus T\gamma_\od$ is indeed a homomorphism.  The inverse homomorphism is $T_1 \oplus T_2 \mapsto T_1+T_2$.
\end{proof}


\begin{lemma}
Under the identifications $K_0 (\scrK(\scrH_+)) \cong \ZZ$ and $K_0 (\scrK(\scrH)^\sigma )  \cong \ZZ\oplus\ZZ$, the map $\average_*$ is  $n \mapsto (n,n)$.
\end{lemma}

\begin{proof}
  Let $p$ be a projection in $\scrK(\scrH_+)$.  Then $p=\chi_{Y_+}p\chi_{Y_+}$, so $p\gamma = \chi_{Y_+}p\gamma\chi_{Y_-}$, and hence $\Tr(p\gamma)=0$.  Since $\gamma_{\ev/\od} = \half(1\pm\gamma)$,  $\Tr(p\gamma_{\ev}) = \Tr(p\gamma_{\od}) = \half \Tr(p)$.  Similarly, $\Tr(\sigma(p)\gamma_{\ev}) =\Tr(\sigma(p)\gamma_{\od}) = \half \Tr(\sigma(p)) = \half \Tr(p)$.  Hence,\linebreak $\Tr(\average(p)\gamma_{\ev}) = \Tr(\average(p)\gamma_{\od}) = \Tr(p)$, and the result follows from the previous lemma.
\end{proof}

By \cite[Proposition 9.4]{Roe:index-theory-book}, $C^*|Y_+|$ has trivial $K$-theory.  The boundary maps in $K$-theory induced from the diagram \eqref{eq:short-exact-sequences} give
\begin{equation}
\label{eqn:boundary-maps}
	\xymatrix{
		K_1(\Roe{Y_+} / \scrK(\scrH_+) ) \ar[r]^\partial_\cong \ar[d]_{\average_*}^\cong &
		K_0(\scrK(\scrH_+)) \ar[d]_{\average_*}  \cong \ZZ &
		n \ar@{|->}[d] 
		\\
		K_1(\Roe{\RR}^\sigma / \scrK(\scrH)^\sigma ) \ar[r]^-\partial &
		K_0(\scrK(\scrH)^\sigma )   \cong \ZZ\oplus\ZZ &
		(n,n).
	}
\end{equation}
We see that $K_1(\Roe{\RR}^\sigma / \scrK(\scrH)^\sigma ) \cong \ZZ$, and the image of its boundary map into $K_0(\scrK(\scrH)^\sigma )$ is $\{(n,n)\st n\in\ZZ\}$.  The corresponding diagram in the other degree gives $K_0(\Roe{\RR}^\sigma / \scrK(\scrH)^\sigma ) \cong 0$.

Now the six-term exact sequence associated to the bottom row of \eqref{eq:short-exact-sequences} becomes
$$
\xymatrix{
	(n,n) & \ZZ\oplus\ZZ \ar[r] &
	K_0( \Roe{\RR}^\sigma ) \ar[r] &
	0 \ar[d] \\
	n \ar@{|->}[u] & \ZZ \ar[u] &
	K_1( \Roe{\RR}^\sigma ) \ar[l] &
	0 \ar[l]
}
$$
Thus, $K_0( \Roe{\RR}^\sigma ) \cong \ZZ$ and $K_1( \Roe{\RR}^\sigma )\cong 0$.  With an appropriate choice of sign, top-left horizontal map is given by $(m,n) \mapsto m-n$.  Applying Lemma \ref {lem:K-theory-of-KHsigma}, this completes the proof of \ref{prop:K-theory}.


\section{The index pairing}

Let $\theta\in K^0(A)$ be the $K$-homology class of a graded unbounded Fredholm module $(H,\rho,D)$, and put $F \defeq D(1+D^2)^{-\half}$.  Let $p$ be a projection in $M_n(A)$.  The index pairing $K^0(A) \times K_0(A) \to \ZZ$ is given by
$$
  ( \theta, [p] ) \defeq \Index \big[ \rho(p)(F\otimes I_n)\rho(p) : \rho(p)H_\ev^n \to \rho(p)H_\od^n \big] ,
$$
(where $I_n$ denotes the identity in $M_n(\CC)$.)

Let $P= \tilde{\rho}(p) \in M_n( \Roe{\RR}^\sigma )$, and let $f$ denote the function $f(x) = x(1+x^2)^{-\half}$, as represented on the geometric $|\RR|$-Hilbert space $\scrH$.  Then
$$
  ( \theta, [p] ) = \Index ( P(f\otimes I_n)P : P\scrH_\ev^n \to P \scrH_\od^n ) .
$$
The right-hand side here depends only on the class of $P$ in $K_0(\Roe{\RR}^\sigma)$.  By Proposition \ref{prop:K-theory}, we may therefore replace $P$ by a finite rank projection $Q$, and the index is
\begin{eqnarray*}
  ( \theta, [p] ) &=& \Index ( Q(f\otimes I_n)Q : Q\scrH_\ev^n \to Q \scrH_\od^n ) \\
  	&=& \dim (Q\scrH_\ev^n) - \dim (Q \scrH_\od^n) \\
	&=& [Q] = \tilde{\rho}_*[p].
\end{eqnarray*}
This completes the proof of Theorem \ref{thm:main-theorem}.

\begin{remark}
Given the above results, it is natural to expect a reformulation of $KK^0(A,\CC)$ in the spirit of Luu.  Indeed, one can define a group $KC^0_\sigma(A,\CC)$ as follows.  Cycles are morphisms from $A$ into the symmetric Roe algebra $\Roe{\RR}^\sigma$.  Equivalence of cycles is generated by controlled unitary equivalences (preserving the involution $\gamma$) and weak homotopies (respecting the symmetry $\sigma$).  Then $KC^0_\sigma(A,\CC) \cong KK^0(A,\CC)$.  We shall not develop this in detail here, as the results follow \cite{Luu:thesis} closely.
\end{remark}


\bibliographystyle{alpha}
\bibliography{dual-control}

\end{document}

%% file: dual-control-macros.tex
\usepackage{amsmath}
\usepackage{amsthm}
\usepackage{amssymb}
\usepackage{verbatim}
\usepackage{ifthen}
\usepackage{amsfonts}
\usepackage[mathscr]{euscript}   
\usepackage{xypic}      
\usepackage{graphicx}
\usepackage{textcomp}
\usepackage{yhmath}
\usepackage{bbm}
\usepackage{empheq}
\usepackage[bookmarks=true, linkbordercolor={0.4 0.35 0.8}]{hyperref}

\newtheorem{theorem}{Theorem}[section]

\newtheorem{lemma}[theorem]{Lemma}
\newtheorem{proposition}[theorem]{Proposition}

\theoremstyle{definition}

\theoremstyle{remark}
\newtheorem{remark}[theorem]{Remark}

\numberwithin{equation}{section}

\newcommand{\scr}{\mathcal}

 \DeclareMathOperator{\Aut}{Aut}

 \DeclareMathOperator{\Tr}{Tr}
 
 \DeclareMathOperator{\Index}{Index}

 \DeclareMathOperator{\Supp}{Supp}

\newcommand{\st}{\; | \;}

\newcommand{\defeq}{:=}

\newcommand{\CC}{\mathbb{C}}

\newcommand{\RR}{\mathbb{R}}

\newcommand{\ZZ}{\mathbb{Z}}

\newcommand{\scrK}{\scr{K}}
\newcommand{\scrB}{\scr{B}}

\newcommand{\half}{\frac{1}{2}}

\newcommand{\pddiff}[3][.]{\ifthenelse{\equal{#1}{.}}
        {\frac{\partial^2 #2}{\partial #3^2}}
        {\frac{\partial^2 #1}{\partial #2 \partial #3}}}

\newcommand{\Univ}[2][.]{\ifthenelse{\equal{#1}{.}}{\mathcal{U}}{\mathcal{U}^{(#1)}}(\mathfrak{#2})}

\newcommand{\ev}{\mathrm{ev}}
\newcommand{\od}{\mathrm{od}}

\newcommand{\Roe}[1]{C^*|#1|}			

\newcommand{\average}{(I+\sigma)}		

\newcommand{\scrH}{\scr{H}}